\documentclass{amsart}[12pt]
\usepackage{amsmath}

\newtheorem{theorem}{Theorem}[section]

\theoremstyle{definition}

\theoremstyle{remark}
\newtheorem{remark}[theorem]{Remark}
\newtheorem*{remarks*}{Remarks}

\numberwithin{equation}{section}

\newcommand{\CC}{{\mathbb C}}

\newcommand{\ZZ}{{\mathbb Z}}
\newcommand{\RR}{{\mathbb R}}
\newcommand{\LL}{{\mathbb L}}

\newcommand{\fg}{{\mathfrak g}}

\newcommand\Cinf{\mathcal{C}^\infty}
\newcommand{\st}[1]{\ensuremath{^{\scriptstyle \textrm{#1}}}}

\newcommand{\hol}{\mbox{\textrm hol}}

\renewcommand{\span}{\mbox{\textrm span}}
\newcommand{\trace}{\mbox{\textrm trace}}

\newcommand{\Vol}{\mbox{\textrm Vol}}

\begin{document}
\openup3pt

\title{The Mellin transform and spectral properties of toric varieties}

\author{Victor Guillemin}\address{Department of Mathematics\\
MIT \\Cambridge, MA 02139 \\ USA}\email{vwg@math.mit.edu}
\author{Zuoqin Wang}
\address{Department of Mathematics\\
MIT \\ Cambridge, MA  02139 \\ USA}\email{wangzq@math.mit.edu}

\begin{abstract}
In this article we apply results of \cite{W} on the twisted Mellin
transform to problems in toric geometry. In particular we use
these results to describe the asymptotics of probability densities
associated with the monomial eigenstates, $z^k$, $k \in \ZZ^d$, in
Bargmann space and prove an ``upstairs" version of the spectral
density theorem of \cite{BGU}. We also obtain for the $z^k$'s,
``upstairs" versions of the results of \cite{STZ} on distribution
laws for eigenstates on toric varieties.
\end{abstract}

\maketitle


\section{Introduction}
\label{sec:1}

Let $X$ be a compact K\"ahler manifold and $\LL \to X$ a
Hermitian line bundle whose curvature form is the negative of the
K\"ahler form.  For every positive integer $N$ let $\LL^N$ be the
$N$\st{th} tensor power of $\LL$ and let $\pi_N$ be orthogonal
projection
\begin{displaymath}
  \Gamma (\LL^N \, ; \, X) \to \Gamma_{\scriptsize \hol}(\LL^N \, ; \, X)\, .
\end{displaymath}
The \emph{spectral measure}, $\mu_{N}$, of the pair $(X,\LL^N)$
is the measure defined by
\begin{equation}
  \label{eq:1.1}
  f \in C (X) \to \mu_N (f) = : \trace\ \pi_N M_f \pi_N
\end{equation}
where $M_f$ is the operator ``multiplication by $f$''.  It is
known that as $N$ tends to infinity this measure has an asymptotic
expansion in inverse powers of $N$ with distributional
coefficients. (See, for instance, \cite{BG}.)  However, the terms
in this expansion have only been computed in a few special cases.
One such case:  $X$ a toric variety, was the topic of a recent
article by Dan Burns, Alejandro Uribe and one of the authors of
this paper, \cite{BGU}.  One of the purposes of this note is to
study the asymptotics of another interesting spectral measure
associated with toric varieties.  More explicitly, by Delzant's
theorem every toric variety $X$ is the GIT quotient of $\CC^d$ by
a subtorus, $G$, of $T^d$, so $\Gamma_{\scriptsize \hol} (\LL^N \,
; \, X)$ has two descriptions: as holomorphic sections of $\LL^N$
and as $G$-invariant holomorphic sections of $(\LL_{\CC^d})^N$
where $\LL_{\CC^d}$ is the line bundle on $\CC^d$ which gives rise
by reduction to the line bundle, $\LL$, on $X$.  In the toric case
this line bundle is the trivial bundle, $\CC^d \times \CC$;
however the Hermitian inner product on it is the (non-trivial)
Bargmann inner product.  Namely for the trivializing section
\begin{displaymath}
  s: \CC^d \to \CC^d \times \CC \, , \quad z \to (z,1)
\end{displaymath}
of $\LL_{\CC^d}$
\begin{equation}
  \label{eq:1.2}
    \langle s,s \rangle = e^{- |z|^2}\, .
\end{equation}

Let $\RR^d$ be the Lie algebra of $T^d$, $\fg$, the Lie algebra of
$G$ and $e_i$, $i=1,\ldots ,d$ the standard basis vectors of
$\RR^d$. From the inclusion, $\fg \to \RR^d$, one gets a dual map
$L: (\RR^d)^* \to \fg^*$ and the weights of the representation of
$G$ on $\CC^d$ are just the vectors
\begin{equation}
  \label{eq:1.3}
  \alpha_i = L e^*_i \, .
\end{equation}
By Delzant's theorem every toric variety is obtained from $\CC^d$
by reducing at a weight, $\alpha \in \ZZ^*_G$, and by the
``quantization commutes with reduction'' theorem
$\Gamma_{\scriptsize \hol} (\LL^N \, ; \, X)$ is isomorphic as
vector space to the space of $L^2$ holomorphic sections of
$(\LL_{\CC^d})^N$ which transform under $G$ according to the
weight $\alpha$, in other words, the space:
\begin{equation}
  \label{eq:1.4}
  \Gamma^N_{\alpha} = \span \{ z^{k_1}_1 \ldots z^{k_d}_d \, , \,
       \sum k_i \alpha_i = N \alpha \} \, .
\end{equation}
Our first goal in this note is to study the spectral measure on
$\CC^d$ defined by
\begin{equation}
  \label{eq:1.5}
  \nu_N (f) = \trace\ \pi_N M_f \pi_N
\end{equation}
where $\pi_N$ is the orthogonal projection of the Bargmann space,
$L^2 (\CC^d \, , \, e^{-N |z|^2} \, dz \, d\bar{z})$, onto
$\Gamma^N_{\alpha}$.  This measure has many features in common with
the measure, $\mu_N$, and in principle, one can be computed from the
other.  However, $\nu_N$ turns out to be much easier to compute.
Moreover, the computation involves an object which is of
considerable interest in itself:  a ``twisted'' version of the
classical Mellin transform.

The asymptotics of (\ref{eq:1.5}) is related via an
Euler-Maclaurin formula to the asymptotics of the measure
   \begin{equation}
   \label{eq:1.6}
   \frac 1{c_{N,k}^2}|z^k|^{2} e^{-N|z|^2}dzd\bar z,
   \end{equation}
where the $c_{N, k}$ is chosen to make this a probability measure;
and the second goal of this paper will be to study this
asymptotics in detail, or, equivalently, to study the asymptotics
of the integral
   \begin{equation}
   \label{eq:1.7}
   \frac 1{c_{N, k}^2} \int |z^k|^{2}e^{-N|z|^2}f dzd\bar z,
   \end{equation}
for $f \in \Cinf_0(\CC^d)$ and $k \in N\Delta$ as $N$ tends to
infinity. The measure (\ref{eq:1.6}) can be thought of as the
probability measure
   \begin{equation*}
   \langle s_k, s_k \rangle dzd\bar z
   \end{equation*}
associated with the ``quantum state"
   \begin{equation*}
   s_k = \frac 1 {c_{N, k}}z^k \in \LL^N_{\CC^d}
   \end{equation*}
in Bargmann space, and our third goal in this paper will be to
study the probability distribution on the real line associated
with this quantum state:
   \begin{equation}
   \label{eq:1.8}
   \sigma_{N, k}([t, \infty)) = \Vol\{z \in \CC^d\ |\ \langle s_k, s_k \rangle(z) \ge
   t\}.
   \end{equation}
The analogue of the distribution (\ref{eq:1.8}) for holomorphic
sections of the line bundle, $\LL$, ``downstairs" on the toric
variety, $X$, was studied in detail by Shiffman-Tate-Zelditch in
\cite{STZ}, and we will show that the distribution (\ref{eq:1.8})
has the same universal rescaling properties as theirs. (In fact
these properties are much more transparently exhibited in this
``upstairs" picture.)

Finally in the last section of this paper we will say a few words
about how these ``upstairs" and ``downstairs" picture are related
via some general results in GIT theory.

\section{The twisted Mellin transform}
\label{sec:2}

For a bounded $\Cinf$ function on $\RR^d$ we define its
``$N$-twisted'' Mellin transform to be the quotient:
\begin{equation}
  \label{eq:2.1}
  A_N f (x) =
     \frac{\int_{\RR^d_+} e^{N (\sum x_i \log y_i -y_i)}f(y) \, dy}
     {\int_{\RR^d_+} e^{N (\sum x_i \log y_i - y_i)} \, dy}\, .
\end{equation}
For properties of this transform and the relation of (\ref{eq:2.1})
to the usual Mellin transform see \cite{W}.  The main property we
will need here is the following asymptotic expansion
\begin{eqnarray}
  \label{eq:2.2}
  A_N f(x) & \sim & \sum_\alpha N^{-|\alpha |}
      f^{(\alpha)}(x)\ g_{\alpha} (Nx),\\
\noalign{\hbox{where}}\nonumber \\
\label{eq:2.3}
g_{\alpha} (x) &=& g_{\alpha_1} (x_1) \ldots g_{\alpha_d}(x_d)\\
\noalign{\hbox{and}}\nonumber \\
\label{eq:2.4}
g_k (s) &=&  \frac 1{k!}\sum_{0 \leq \ell \leq k}
(-1)^{\ell} \binom{k}{\ell}
   s^{\ell} s^{(k-\ell)},
\end{eqnarray}
$s^{(m)}$ being the Hardy function
\begin{equation}
  \label{eq:2.5}
s^{(m)} = (s+m) \ldots (s+1) \, .
\end{equation}
The sequence of functions $g_k(s)$ also has a nice generating
function description in term of identity
\begin{equation}
\label{gf} \sum_r g_k(s)x^k = \frac{e^{-sx}}{(1-x)^{1+s}}.
\end{equation}
That (\ref{eq:2.2}) \emph{is} in fact an asymptotic expansion
follows from one of the main results of \cite{W}.

\begin{theorem}
  \label{th:2.1}
The $k$\st{th} degree polynomial (\ref{eq:2.4}) is in fact a
polynomial of degree $[\frac{k}{2}]$.
\end{theorem}

Thus the $\alpha$\st{th} summand in (\ref{eq:2.2}) is of order
$N^{-\frac{|\alpha |}{2}}$.

We will now show how this asymptotics is related to the asymptotics
of the measure~(\ref{eq:1.5}).  Let $\Delta$ be the convex polytope
\begin{equation}
  \label{eq:2.6}
  \{ x \in \RR^d_+ \, , \, \sum x_i \alpha_i = \alpha \}\, .
\end{equation}
Then the functions, $z^k$, $k \in N \Delta \cap \ZZ^d$, are an
orthogonal basis of $\Gamma^N_{\alpha}$ and the functions
\begin{equation}
  \label{eq:2.7}
  \frac{1}{c_{N,k}} z^k
\end{equation}
with
\begin{equation}
  \label{eq:2.8}
  c_{N,k} = \left( \int_{\CC^d} |z^{k}|^2 e^{-N |z|^2}\, dz \,
    d\bar{z} \right)^{\frac{1}{2}}
\end{equation}
are an orthogonal basis of $\Gamma^N_{\alpha}$.  Hence the trace
of $\pi_N M_f \pi_N$ is equal to the sum over $k \in N \Delta
\cap \ZZ^d$ of
\begin{equation}
  \label{eq:2.9}
  \frac{\int |z_1|^{2k_1} \ldots |z_d|^{2k_d} e^{-N|z|^2}
             f(z) \, dz \, d\bar{z}}
        { \int |z_1|^{2k_1} \ldots | z_d|^{2k_d} e^{-N|z|^2}\,
          dz \, d\bar{z}}\, .
\end{equation}
Now note that  $\nu_N$ is $T^d$-invariant, so to compute $\nu_N (f)$
it suffices to compute $\nu_N (f)$ for functions which are
themselves $T^d$-invariant, i.e.,~functions of the form $f (|z_1|^2
, \ldots, |z_d|^2)$  where $f (x_1,\ldots ,x_n)$ is a bounded
$\Cinf$ function on $\RR^d$.  However for such functions,
(\ref{eq:2.9}) becomes
\begin{equation}
  \label{eq:2.10}
  \frac{\int_{\RR^d_+} x^{k_1}_1 \ldots x^{k_d}_d
       e^{-N (x_1 + \cdots + x_d)} f(x) \, dx}
   {\int_{\RR^d_+} x^{k_1}_1 \ldots x^{k_d}_d
       e^{-N (x_1 + \cdots + x_d)} \, dx} \, .
\end{equation}
This shows that the integral (1.7) is just the twisted Mellin
transform, $A_N(f)(x)$ evaluated at $x = \frac kN$, and hence
gives us for this integral the asymptotic expansion (2.2).
Moreover, by summing (1.7) over the lattice points, $k \in N\Delta
\cap \ZZ^d$ we get for the spectral measure, $\nu_N$, the formula
\begin{equation}
  \label{eq:2.11}
  \nu_N (f) = \sum_{k \in N \Delta \cap \ZZ^d} A_N f
     \left( \frac{k}{N}\right) \, .
\end{equation}

In the next section we will get an asymptotic expansion for
$\nu_N (f)$ by interpreting the sum on the right as a Riemann sum
and combining the formula~(\ref{eq:2.2}) with an Euler--Maclaurin
formula for Riemann sums over convex polytopes.

\section{Riemann sums}
\label{sec:3}

Let $\Delta \subseteq \RR^n$ be an $n$-dimensional convex
polytope.  By elementary calculus the Riemann integral of a
function $f \in \Cinf (\Delta)$:
\begin{displaymath}
  \int_{\Delta} f(x) \, dx
\end{displaymath}
is approximated by the Riemann sum
\begin{displaymath}
  \frac{1}{N^n} \sum_{k \in N \Delta \cap \ZZ^n} f \left( \frac{k}{N}\right) \,
\end{displaymath}
up to an error term of order $O (N^{-1})$.

Recently Guillemin and Sternberg showed that if $\Delta$ is a
\emph{lattice} polytope, i.e.,~if its vertices are lattice points,
then this $O (N^{-1})$ can be replaced by an asymptotic series in
inverse powers of~$N$.    In particular for polytopes associated
with toric varieties (such as the polytope (\ref{eq:2.6})) the terms
in this series can be explicitly computed by the following method.

Enumerate the facets of $\Delta$, and for the $i$\st{th} facet
let $u_i \in \ZZ^n$ be a primitive lattice vector which is
perpendicular to this facet and points ``outward'' from $\Delta$
into $\RR^n$.  Then $\Delta$ can be defined by a set of
inequalities
\begin{equation}
  \label{eq:3.1}
  \langle u_i , x\rangle \leq c_i \, , \quad i=1,\ldots ,r
\end{equation}
where $r$ is the number of facets.  Let $\Delta_h$ be the polytope
\begin{equation}
  \label{eq:3.2}
  \langle u_i , x\rangle \leq c_i  + h_i \, ,
       \quad i=1,\ldots ,r \, .
\end{equation}
Then for $f \in \Cinf (\RR^n)$
\begin{eqnarray}
  \label{eq:3.3}
  \frac{1}{N^n} \sum_{k \in \ZZ^n \cap N \Delta} f
      \left(\frac{k}{N}\right) & \sim &
         \left( \tau \left(\frac{1}{N}\, \frac{\partial}{\partial h}\right)
         \int_{\Delta_h} f(x) \, dx \right) (h=0)\\
\noalign{\hbox{where}}\nonumber\\
\label{eq:3.4}
 \tau (w_1 , \ldots, w_r ) & = & \prod^r_{i=1}
      \frac{w_i}{1-e^{-w_i}}
\end{eqnarray}
and $\tau \left( \frac{1}{N}\, \frac{\partial}{\partial
    h}\right)$ is the operator obtained from (\ref{eq:3.4}) by
making the substitution $w_i \to \frac{1}{N}\, \frac{\partial}{\partial
    h_i}$.  (The leading term on the right hand side of
  (\ref{eq:3.3}) is the integral of $f$ over $\Delta$ and the
  higher order terms are the asymptotic series in inverse powers
  of $N$ that we just alluded to.)

Now notice that if we divide (\ref{eq:2.11}) by $N^n$ the
  right hand side is exactly a Riemann sum of the form above.
  Hence if we replace $A_N f$ by the series (\ref{eq:2.2})
  and apply (\ref{eq:3.3}) to each summand we get an asymptotic
  expansion of $\nu_N (f)$ in inverse powers of $N$ in which the
  summands can be read off from the summands on the right hand
  side of (\ref{eq:3.3}).  (For a more detailed description of
  these summands see \S3 of \cite{GS06}).

\section{The Nonrescaled Distribution Law}

We turn now to the third topic of this paper: the asymptotics of
the probability distribution (1.8). Suppose $k=Na$ with $a \in
\Delta$. We begin by observing that
   \begin{equation*}
   \|z^{k}\|^2_N = \int_{\CC^d}\langle z^{k}, z^{k} \rangle_N dzd\bar z
    = \left( \frac{\pi}N \right)^d \prod_i \frac{(k_i)!}{N^{k_i}},
   \end{equation*}
and hence
   \begin{equation}
   \label{eq:4.1}
   \langle s_{k}, s_{k} \rangle_N =
   \left(\frac{N}{\pi} \right)^d \frac{N^{|k|}}{k!}
   |z^k|^{2} e^{-N|z|^2}.
   \end{equation}

We first assume that $k=(k_1, \cdots, k_d)$ with $k_i > 0$ for all
$1 \le i \le d$, and observe that $\sigma_{N,k}([t, \infty))$ is
the volume of the region in $\CC^d$
   \begin{equation}
   \label{eq:4.2}
   |z^k|^{2} e^{-N|z|^2} > \left(\frac {\pi} N\right)^d \frac{k!}{N^{|k|}}t,
   \end{equation}
or, with $a=\frac kN$, the region
   \begin{equation}
   \label{eq:4.3}
   |z^a|^{2} e^{-|z|^2} > \left( \left(\frac {\pi}N\right)^d \frac{k!}{N^{|k|}}
   t\right)^{1/N}.
   \end{equation}

By Stirling's formula,
   \begin{equation*}
   k_i! = \sqrt{2 \pi k_i}\left( \frac{k_i}e \right)^{k_i}\left(1+O(\frac 1N)\right),
   \end{equation*}
so the right hand side of (\ref{eq:4.3}) is equal to
   \begin{equation*}
   \lambda_N = \left(\pi^d  N^{-d/2} t \prod_i(2\pi a_i)^{1/2} \right)^{1/N} \left(
   \frac{a}{e}\right)^a \left(1+O(\frac 1{N^2})\right).
   \end{equation*}
Thus if we set $|z_i|^2=r_i$ and let $f(r)$ be the function
   \begin{equation*}
   f(r)=\sum_{i=1}^d \left(a_i \log r_i - r_i\right),
   \end{equation*}
the inequality (\ref{eq:4.3}) becomes
   \begin{equation}
   \label{eq:4.4}
   f(r) \ge \log{\lambda_N} = \sum_{i}\left(a_i \log a_i - a_i\right) - \frac d{2N}\log{N}
   +\frac{\log{t}+\gamma}{N} +
   O(\frac 1{N^2}),
   \end{equation}
where
   \begin{equation}
   \label{eq:4.5}
   \gamma =  \log\left(\pi^d \prod_i (2\pi a_i)^{1/2}\right).
   \end{equation}
We now note that $f(r)$ has a unique maximum at $r=a$ and that in
a neighborhood of this maximum,
   \begin{equation*}
   f(r) = \sum_i \left(a_i \log a_i -a_i - \frac 1{2a_i}(r_i-a_i)^2\right) + \cdots.
   \end{equation*}
Hence for $N$ large (ignoring terms in $N$ of order $O(\frac 1N)$)
(\ref{eq:4.3}) reduces to
   \begin{equation*}
   \sum_{i} \frac 1{2a_i}(r_i-a_i)^2 \le \frac d{2N}\log N,
   \end{equation*}
or, since $r_i = |z_i|^2$,
   \begin{equation}
   \label{eq:4.6}
   \sum_{i} \frac 1{2a_i} \left(|z_i|^2 - a_i\right)^2 \le \frac
   d{2N} \log{N} + O(\frac 1N).
   \end{equation}
To compute the volume of this set to the leading order, we first
note that the volume of the ellipsoid
   \begin{equation}
   \label{eq:4.7}
   \sum_{i=1}^d \frac 1{2a_i} x_i^2 \le \varepsilon
   \end{equation}
in $\RR^d$ is
   \begin{equation}
   \label{eq:4.8}
   \gamma_d \left(\prod 2a_i\right)^{1/2} \varepsilon^{d/2},
   \end{equation}
where $\gamma_d$ is the volume of the unit $d$-ball. Now consider
the map
   \begin{equation*}
   g: \RR_+^d \to \RR^d, \qquad s_i \mapsto x_i=s_i^2 -a_i.
   \end{equation*}
The pre-image of the region (\ref{eq:4.7}) with respect to this
map is the set
   \begin{equation}
   \label{eq:4.9}
   \sum \frac 1{2a_i}(s_i^2-a_i)^2 \le \varepsilon.
   \end{equation}
If $s$ is a point in this set, then $s_i = \sqrt{a_i} +
O(\varepsilon^{1/4})$, so
   \begin{equation*}
   \det(Dg(s)) = \prod_i (2s_i) = 2^d \prod_i \sqrt{a_i},
   \end{equation*}
and thus by (\ref{eq:4.8}) the volume of the region (\ref{eq:4.9})
is equal, modulo $O(\varepsilon^{1/4})$, to
   \begin{equation}
   \label{eq:4.10}
   \gamma_d \left(\frac \varepsilon 2 \right)^{d/2}.
   \end{equation}

Finally note that the region (\ref{eq:4.6}) is, with $\varepsilon
= \frac d{2N}\log{N}$, the pre-image of the region (\ref{eq:4.9})
with respect to the torus fibration, $s_i=|z_i|$. Since each torus
fiber has volume $\prod (2\pi s_i)$ and $s_i =
\sqrt{a_i}+O(\varepsilon^{1/4})$, the total volume of the region
(\ref{eq:4.6}) is equal modulo a factor of
$1+O(\varepsilon^{1/4})$ to
   \begin{equation}
   \label{eq:4.11}
   (2\pi)^d \gamma_d \left(\prod_i
   \frac{a_i\varepsilon}2\right)^{1/2},
   \end{equation}
and hence by substituting $\frac{d \log{N}}{2N}$ for $\varepsilon$
we arrive finally at the asymptotic formula
   \begin{equation}
   \label{eq:4.12}
   \sigma_{N, k}\left([t, \infty)\right) \sim \pi^d \gamma_d
   \prod_i \left(a_i \frac dN \log{N}\right)^{1/2}.
   \end{equation}

\begin{remark}
More generally suppose $k=(k_1, \cdots, k_l, 0, \cdots, 0)$ with
$k_i>0$ for $1 \le i \le l$, then $\sigma_{N, k}([t, \infty))$ is,
to its leading order, equal to the volume of the region
   \begin{equation*}
   \sum_{i=1}^l \frac 1{2a_i}(|z_i|^2-a_i)^2 + \sum_{i=l+1}^d
   |z_i|^2 \le \frac d{2N}\log{N}
   \end{equation*}
To compute the volume of this set, we  regard it as the pre-image
of the $l$-torus fibration over the $2d-l$ dimensional ellipsoid
   \begin{equation*}
   \sum_{i=1}^l \frac 1{2a_i}(s_i^2-a_i)^2 + \sum_{i=l+1}^d
   (x_i^2+y_i^2)\le \frac d{2N}\log{N},
   \end{equation*}
and by the same argument as above, get
   \begin{equation}
   \label{eq:4.13}
   \sigma_{N, k}\left([t, \infty)\right) \sim 2^{l-d} \pi^l
   \gamma_{2d-l} \left( \frac {d\log{N}}N\right)^{d-\frac l2}
   \prod_i \left(a_i\right)^{1/2}.
   \end{equation}
\end{remark}

\section{Rescaled Distribution Laws}

For simplicity we assume all $k_i$ are positive. From
(\ref{eq:4.3}), (\ref{eq:4.4}) and (\ref{eq:4.5}) we have
   \begin{equation}
   \label{eq:5.1}
   \sum_{i=1}^d \frac 1{2a_i} \left(|z_i|^2 - a_i\right)^2 \le
   \varepsilon_N,
   \end{equation}
where
   \begin{equation}
   \label{eq:5.2}
   \varepsilon_N = \frac d{2N} \log N - \frac{\log{t}+\gamma}N +
   o\left(\frac 1{N}\right).
   \end{equation}
Thus the $t$ term gets absorbed in the $O(\frac 1N)$ and doesn't
affect the leading asymptotics of $\sigma_N([t, \infty))$.
However, we can remedy this problem by rescaling techniques.

The first choice of rescaling is to eliminate the leading term
$\frac{d}{2N}\log{N}$. To do so, we replace $t$ by $N^{d/2}t$.
Then
   \begin{equation}
   \label{eq:5.3}
   \varepsilon_N = (-\log{t}-\gamma)\frac 1N + o(\frac 1N)
   \end{equation}
and the computations in the last section show that this rescaled
version of $\sigma_{N,k}([t, \infty))$ satisfies (\ref{eq:4.9})
with $\varepsilon = \varepsilon_N$ given by (\ref{eq:5.4}) and
hence depends in an interesting way on $t$. (One proviso, however, is that
$\log{t}$ has to be smaller than $-\gamma$.)

There are also many other interesting choices of rescalings: we
may rescale $t$ such that the term containing $\log{t}$ dominate
other terms. For example, we may replace $t$ by $e^{-N^\alpha
(\log{N})^\beta t}$, where $0< \alpha < 1$ or $\alpha=0, \beta >
1$. In this case
   \begin{equation}
   \label{eq:5.4}
   \varepsilon_N = N^{\alpha-1}(\log{N})^\beta t +
   O(\frac{\log{N}}N).
   \end{equation}
We may also replace $t$ by $N^{-t}$, which is the extreme case
$\alpha=0, \beta=1$ above, then
   \begin{equation}
   \label{eq:5.6}
   \varepsilon_N = \frac{d+2t}N \log{N} + O(\frac{1}N).
   \end{equation}

\section{``Upstairs" versus ``downstairs"}

In GIT theory one is given a projective variety, $M$, a positive
line bundle, $\LL \to M$, and an action of an algebraic group,
$G_\CC$, on the pair $(M, \LL)$; and one's goal is to make sense
of the quotients, $M/G_\CC$ and $\LL/G_\CC$. Unfortunately
$M/G_\CC$ can be a fairly pathological object, even as a
topological space (e.g. it can be non-Hausdorff); however, $M$
contains a $G_\CC$ invariant Zariski open set, $M_{stable}$ for which the
quotient
   \begin{equation}
   \label{eq:6.1}
   M_{red} = M_{stable}/G_\CC
   \end{equation}
is a quasi-projective variety. Moreover if $\LL_{stable}$ is the
restriction of $\LL$ to $M_{stable}$ one gets a quotient line
bundle
   \begin{equation}
   \label{eq:6.2}
   \LL_{red} = \LL_{stable}/G_\CC
   \end{equation}
on $M_{red}$, and modulo some Kaehlerizability assumptions one can
equip $\LL$ and $\LL_{red}$ with natural Hermitian inner products,
$\langle \cdot, \cdot \rangle$ and $\langle \cdot, \cdot
\rangle_{red}$ which are related by an identity
   \begin{equation*}
   \pi^* \langle\cdot, \cdot \rangle_{red} = e^\psi \langle \cdot,
   \cdot \rangle
   \end{equation*}
where $\pi$ is the projection of $M_{stable}$ onto $M_{red}$ and
$\psi$ is a $\Cinf$ function called the stability function. (See
\cite{GS82}, \cite{BGU}, \cite{HK}, \cite{L} etc.)

In the case of toric varieties $M$ is $\CC P^{d-1}$ (or more
simply its de-projectivization, $\CC^d$), $G_\CC$ is a complex
torus and the action of $G_\CC$ on $\CC^d$ is a linear action; and
in the toric variety case the stability function turns out to be
given by a fairly simple formula: If $X=M_{red}$ is the reduction
of $\CC^d$ at a weight, $\alpha$, of $G$, $\psi$ is of the form
  \begin{equation}
  \label{eq:6.3}
  \psi = -|z|^2 + f,
  \end{equation}
where $f$ transforms under the action of $G_\RR$ according to the
character, $e^{\alpha}$. From this description of $\psi$ one can
conclude quite a bit about the asymptotic behavior of norm-squares
of holomorphic sections of $\LL^N_X$ from the ``upstairs picture"
on $\CC^d$ that we've described in this paper. In particular one
can obtain some of the results of Shiffman-Tate-Zelditch that we
mentioned in the introduction by combining results about ``moments
of measures" proved by them in $\S 4.1$ of their paper with
asymptotic properties of $e^{N\psi}$ and results of $\S 5$ above.
We will give an account of this stability theory for toric
varieties in \cite{GW}.

\end{document}